# Coplanar Low Thrust Transfer with Eclipses Using Analytical Costate Guess

Max CERF [*]
Ariane Group

## Abstract

Low thrust orbital transfers are difficult to optimize by indirect methods. The main issues come from the costate guess and from the numerical propagation accuracy required by the shooting method. In the case of a coplanar minimum-time low thrust transfer with eclipses, an analytical costate guess is proposed. The optimal control problem reduces to an unconstrained minimization problem with two unknowns. A derivative free algorithm yields a quasi-optimal solution from scratch in a few minutes. No specific guess is necessary and the algorithm properties ensure finding the global minimum of the unconstrained problem. The method is applicable to thrust levels from large to very low and to any eclipse configuration, as exemplified on a transfer towards the geostationary orbit.

## 1. Introduction

Electric rocket engines offer promising propellant savings for orbital transfers at the expense of long durations due to very low thrust levels. Minimum-time low thrust trajectories are consequently studied intensively for several decades. An additional difficulty comes from the large power required by electric propulsion engines, impeding their use during the transit in the Earth shadow. The minimum-time low thrust transfer between two orbits formulates as an optimal control problem with dynamics discontinuities at the eclipse entrances and exits. Various approaches have been applied to solve such optimal control problems as efficiently as possible. They are classified between direct and indirect methods [1-3].

Direct methods discretize the optimal control problem in order to rewrite it as a nonlinear large scale optimization problem. Various discretization methods can be chosen for the dynamics and for the control [4-9]. This process is straightforward and new variables or constraints may be added to the problem with reduced programming effort.

---
[*] Ariane Group, 78130 Les Mureaux, France
max.cerf@ariane.group

Several software packages such as IPOPT, BOCOP, GESOP, SNOPT, WORHP …, are available to solve the large scale optimization problem. The direct approach is suitable to a wide range of applications. The main drawback is that it is computationally expensive when the number of variables becomes large, which is especially the case for low-thrust transfers [10-13]. Finding an accurate solution for such problems may be difficult. An alternate discretization approach which reduces the problem size consists in approximating the dynamics by a series of impulsive maneuvers [13-16].

On the other hand indirect methods are based on the Pontryagin Maximum Principle (PMP) [17-18] which gives a set of necessary conditions for a local optimal solution. For minimum-time problems the PMP yields the optimal thrust direction aligned with the velocity costate [19-22]. The problem is reduced to a Two Point Boundary Value Problem (TPBVP). The resulting nonlinear system is generally solved by a shooting method using a Newton-like algorithm. The convergence is fast and accurate, but the method requires both an adequate starting point and a high integration accuracy. These are major issues when applying the indirect approach to low thrust transfers [23,24], particularly in the case of dynamics discontinuities [25]. Moreover the problem may admit singular solutions for which the PMP first order conditions no longer define the optimal control. Such singular solutions require further theoretical analysis and specific solution methods [26,27].

For low thrust transfers various approaches can be envisioned to build a satisfying initial costate guess and benefit from the efficiency of the indirect method. In [28] the impulse transfer solution is used to provide a good initial guess to the shooting algorithm. This method is based on the fact that a continuous high-thrust orbit transfer shares similarities with the impulse transfer as outlined in [29,30]. Analytical costate approximation are derived in [31] for transfers between circular orbits. In [32] the similarity between the double integrator and the orbital transfer is exploited to propose an analytical costate guess. Multiple shooting reduces the overall sensitivity by splitting the trajectory in several arcs at the expense of additional unknowns and boundary conditions. In [33] a multiple shooting method parameterized by the number of thrust arcs is used to solve an Earth-Mars transfer. The multiple shooting is combined with a collocation method in [34] by splitting the trajectory into thrust and coast arcs in the Earth shadow. Homotopic approaches [35] solve a series of optimization problems by continuous transformation starting from a known solution. In [23,36] a differential continuation method linking the minimization of the $L^2$-norm of the control to the minimization of the consumption is used to solve the low-thrust orbit transfer around the Earth. In [37] simplified formulas are established by interpolating many numerical experiments, which allows a successful initialization for the minimal time orbit transfer problem, in a given range of nearly circular initial and final orbits. Based on that initial guess and on averaging techniques, the software T3D [38] implements continuation and smoothing processes in order to solve minimal time or minimal fuel consumption orbit transfer problems. Particle swarm [39], genetic algorithms [40] or other metaheuristics [41] can also be used to explore largely the variables space and produce a satisfying initial solution. We can also mention mixed methods that use a discretization of the PMP necessary conditions and then apply a large-scale equation solver [42] and dynamic programming methods that search for the global optimum in a discretized state space by solving the Hamilton-Jacobi-Bellman equation [43].

This paper addresses the minimum-time orbital transfer in the coplanar case. The engine is assumed thrusting at a constant level excepted during eclipses. The optimal control problem is addressed by an indirect method using the optimal thrust direction derived from the PMP. The costate guess issue is overcome using an analytical form derived from previous results obtained for high thrust transfers. The unknowns reduce to two angles defining the initial costate vector. The minimum-time problem is then reformulated as an unconstrained minimization problem whose solution can be found from scratch by a derivative free algorithm. By this way an accurate numerical propagation is no longer mandatory as is the case with Newton-like methods. This allows using rectangular coordinates (which makes the dynamics equations and the eclipse conditions simple), large time steps and less stringent tolerance on the numerical integration error.

The text is organized as follows. Section §2 formulates the optimal control problem based on a dynamical model with eclipses. The extremal conditions are analyzed to derive explicitly the costate discontinuities at the eclipse bounds. Section §3 presents the solution method with the analytical costate guess, the unconstrained problem formulation and the features of the DIRECT optimization algorithm. Section §4 presents an application for a typical transfer towards the geostationary orbit considering various eclipse conditions and various thrust levels.

## 2. Problem Formulation and Analysis

This section describes the dynamical model for the planar low thrust transfer with eclipses. The Optimal Control Problem (OCP) is then analyzed by applying the Pontryaguin Maximum Principle (PMP).

### 2.1 Dynamical model

The problem consists in finding the minimum time trajectory to transfer the vehicle from an initial orbit to a final orbit using a low thrust engine which must be switched off during eclipses. The initial and final orbits are coplanar and defined by their apogee and perigee altitudes. The Earth is modeled as a sphere, with its center at the origin of an inertial frame. The Sun is modeled as a point rotating around the Earth at a constant rotation rate. The Sun attraction is neglected and only the Sun direction is used to determine the Earth shadow region.

The vehicle is modeled as a material weighting point with position $\vec{r}(t)$, velocity $\vec{v}(t)$, mass m(t) submitted to the Earth acceleration gravity denoted $\vec{g}(\vec{r})$ and to the engine thrust. The thrust direction can be chosen freely and it is orientated along the unit vector $\vec{u}(t)$. The thrust level T is constant with a burned propellant exhaust velocity equal to $v_e$. The engine is ignited at the initial date $t_0$ and it is constantly thrusting until the final date $t_f$ except during eclipses where it must be switched off.

The eclipse conditions depend on the respective positions of the vehicle $\vec{r}(t)$ and the Sun $\vec{r}_S(t)$ with respect to the Earth as depicted on the Figure 1.

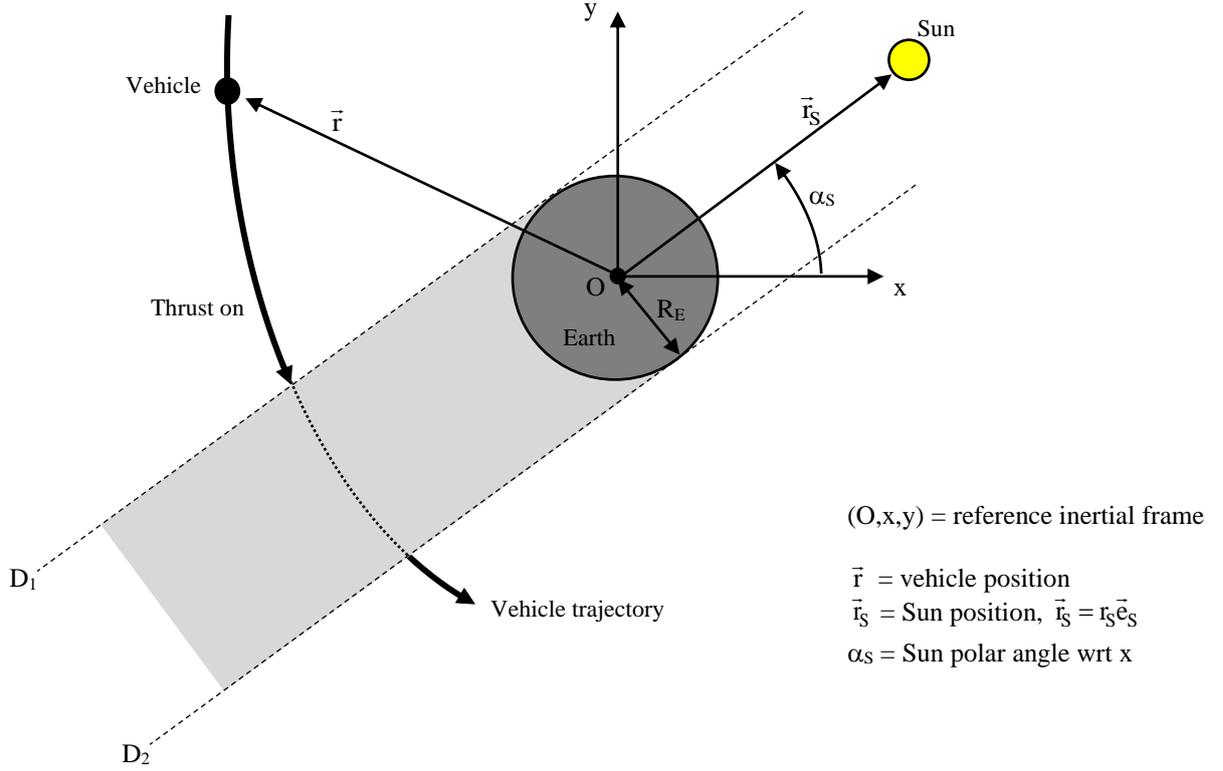

Figure 1 : Eclipse conditions

The Earth shadow is assumed to be a cylinder with the radius of the Earth $R_E$. An alternative model is a conic shape centered at the Sun. For the orbital transfers under consideration, the difference between these shadow models is not significant.

Denoting $\vec{e}_S$ the unit vector pointing toward the Sun, the eclipse conditions under vector form are

$$\begin{cases} \|\vec{e}_S \wedge \vec{r}\| < R_E \\ \vec{e}_S \cdot \vec{r} < 0 \end{cases} \quad (1)$$

The first condition express that the vehicle is inside the shadow cylinder, the second condition checks that the vehicle and the Sun are in opposite directions wrt the Earth.

In this paper, we limit ourselves to planar transfers. The motion equations are written in an inertial reference frame (O,x,y) with the origin at the Earth center. In order to express the eclipse conditions with the coordinates (x,y), we denote $D_1$ and $D_2$ the straight lines limiting the shadow region, $\alpha_S$ the Sun polar angle wrt the x axis and $\omega_S$ the constant angular of the Sun motion in the Earth inertial frame (O,x,y). Starting from the position $\alpha_0$ at the initial date $t_0=0$, the Sun motion is defined by

$$\alpha_S(t) = \alpha_0 + \omega_S t \quad (2)$$

The equations of the lines $D_1$ and $D_2$ are respectively

$$\begin{aligned} D_1: & \quad d_1(x,y) = x\sin\alpha_S - y\cos\alpha_S + R_E = 0 \\ D_2: & \quad d_2(x,y) = x\sin\alpha_S - y\cos\alpha_S - R_E = 0 \end{aligned} \quad (3)$$

The sign of the functions $d_1(x,y)$ and $d_2(x,y)$ determine the position relatively to the shadow region bounds. An eclipse occurs if $d_1(x,y)$ and $d_2(x,y)$ have opposite signs and the position vectors $\vec{r}(t)$ and $\vec{r}_S(t)$ make an obtuse angle. Consistently with Eq. (1) the eclipse conditions in Cartesian coordinates are

$$\begin{cases} d_1(x,y)d_2(x,y) < 0 \\ \vec{r}.\vec{r}_S < 0 \end{cases} \Leftrightarrow \begin{cases} (x\sin\alpha_S - y\cos\alpha_S)^2 - R_E^2 < 0 \\ x\cos\alpha_S + y\sin\alpha_S < 0 \end{cases} \tag{4}$$

When the vehicle enters the Earth shadow, the engine is switched off. Applying the fundamental dynamics principle in the Earth-centered inertial frame yields the motion equations respectively in the light and shadow regions. The dependencies on time (for $\vec{r}$, $\vec{v}$, m and $\vec{u}$) are omitted for conciseness.

$$\begin{cases} \dot{\vec{r}} = \vec{v} \\ \dot{\vec{v}} = \vec{g}(\vec{r}) + \dfrac{T}{m}\vec{u} \quad \text{in light} \\ \dot{m} = -\dfrac{T}{v_e} \end{cases} \quad \text{and} \quad \begin{cases} \dot{\vec{r}} = \vec{v} \\ \dot{\vec{v}} = \vec{g}(\vec{r}) \quad \text{in shadow} \\ \dot{m} = 0 \end{cases} \tag{5}$$

The dynamics is discontinuous at the entrance and exit of the shadow region. The discontinuities occur whenever the function $d_1(x,y)$ or the function $d_2(x,y)$ vanishes, with the additional condition on the vehicle and the Sun being in opposite directions. At a discontinuity date denoted $t_d$, the vehicle position $(x(t_d),y(t_d))$ and the Sun direction $\alpha(t_d)$ satisfy the constraint

$$\psi_S[x(t_d),y(t_d),t_d] \pm R_E = 0 \tag{6}$$

with the function $\psi_S$ defined by

$$\psi_S(x,y,t) = x\sin(\alpha_0 + \omega_S t) - y\cos(\alpha_0 + \omega_S t) \tag{7}$$

This condition will be treated as an interior point constraint with dynamics discontinuity for the optimal control problem.

### 2.2 Optimal control problem

The Optimal Control Problem (OCP) is formulated considering as state variables $\vec{r}(t)$, $\vec{v}(t)$, m(t) and as control variables the thrust direction $\vec{u}(t)$ and the final time $t_f$. The initial state is completely prescribed. The final state is constrained by the targeted apogee and perigee altitudes denoted respectively $h_A$ and $h_P$. The apogee and perigee altitudes actually achieved at the final date are denoted respectively $\psi_A$ and $\psi_P$ and they depend on the final position $\vec{r}(t_f)$ and velocity $\vec{v}(t_f)$. The cost is the final time to minimize.

The OCP formulation is

$$\min_{\vec{u}(t),t_f} J = t_f \quad \text{s.t.} \quad \begin{cases} \dot{\vec{r}} = \vec{v} \\ \dot{\vec{v}} = \vec{g} + \dfrac{\varepsilon T}{m}\vec{u} \\ \dot{m} = -\dfrac{\varepsilon T}{v_e} \end{cases} \quad \text{with} \quad \begin{cases} \vec{r}(t_0) = \vec{r}_0 \\ \vec{v}(t_0) = \vec{v}_0 \quad \text{fixed initial state} \\ m(t_0) = m_0 \end{cases} \\ \begin{cases} \psi_A[\vec{r}(t_f),\vec{v}(t_f)] = h_A \\ \psi_P[\vec{r}(t_f),\vec{v}(t_f)] = h_P \end{cases} \quad \text{constrained final state} \tag{8}$$

The variable ε is equal to 0 (respectively 1) inside (respectively outside) the shadow region. The dynamics is discontinuous at the shadow entrance and exit defined by the interior point constraints Eq. (6). The number of such discontinuities is unknown a priori since it depends on the number of revolutions performed during the transfer. It can be noticed that the problem is not autonomous due to the Sun position $\vec{r}_S$ which depends explicitly on the time. The optimal solution depends thus on the initial Sun position defined by the angle $\alpha_0$.

### 2.3 Extremal Analysis

The optimal trajectory is sought by applying the Pontryaguin Maximum Principle (PMP) [17,18]. For that purpose we introduce the costate vectors $\vec{p}_r(t)$, $\vec{p}_v(t)$, $p_m(t)$ associated respectively to the position, the velocity, the mass. These costate vectors do not vanish identically on any interval of $[t_0, t_f]$ and they are defined up to a non-positive scalar multiplier $p_0$. We choose the usual normalization for regular extremals $p_0 = -1$.

With these notations, the Hamiltonian for the OCP Eq. (8) is

$$H = \vec{p}_r \cdot \dot{\vec{r}} + \vec{p}_v \cdot \dot{\vec{v}} + p_m \dot{m} = \vec{p}_r \cdot \vec{v} + \vec{p}_v \cdot \left( \vec{g} + \frac{\varepsilon T}{m} \vec{u} \right) + p_m \left( -\frac{\varepsilon T}{v_e} \right)$$
$$= \vec{p}_r \cdot \vec{v} + \vec{p}_v \cdot \vec{g} + \varepsilon T \left( \frac{\vec{p}_v \cdot \vec{u}}{m} - \frac{p_m}{v_e} \right) \quad (9)$$

The PMP provides the following first order necessary conditions on $\vec{u}$ and $t_f$ to be an optimal control.

- The Hamiltonian maximization condition with respect to the control $\vec{u}(t)$.

$$\max_{\vec{u}} H \quad (10)$$

This condition leads to a thrust direction aligned with the velocity costate.

$$\vec{u} = \frac{\vec{p}_v}{p_v} \quad \text{with} \quad p_v = \|\vec{p}_v\| \quad (11)$$

and the Hamiltonian Eq. (9) can be written as

$$H = \vec{p}_r \cdot \vec{v} + \vec{p}_v \cdot \vec{g} + \varepsilon T \Phi \quad \text{with} \quad \Phi = \frac{p_v}{m} - \frac{p_m}{v_e} \quad (12)$$

- The costate differential equations.

$$\begin{cases} \dot{\vec{p}}_r = -\frac{\partial H}{\partial \vec{r}} = -\frac{\partial \vec{g}}{\partial \vec{r}} \vec{p}_v \\ \dot{\vec{p}}_v = -\frac{\partial H}{\partial \vec{v}} = -\vec{p}_r \\ \dot{p}_m = -\frac{\partial H}{\partial m} = \frac{\varepsilon T}{m^2} p_v \end{cases} \quad (13)$$

- The transversality conditions on the final costate, derived from the final constraints $\psi_A$, $\psi_P$ with the respective multipliers $\nu_A$, $\nu_P$ and from the final cost $t_f$.

$$\begin{cases} \vec{p}_r(t_f) &= \nu_A \dfrac{d\psi_A}{d\vec{r}}(t_f) + \nu_P \dfrac{d\psi_P}{d\vec{r}}(t_f) \\ \vec{p}_v(t_f) &= \nu_A \dfrac{d\psi_A}{d\vec{v}}(t_f) + \nu_P \dfrac{d\psi_P}{d\vec{v}}(t_f) \\ p_m(t_f) &= 0 \end{cases} \qquad (14)$$

- The transversality condition on the free final date $t_f$.

$$H(t_f) = -1 \qquad (15)$$

The dynamics discontinuity at the shadow entrance or exit generates discontinuities on the costate components and on the Hamiltonian [44,45]. A discontinuity occurs whenever the condition Eq. (6) is met with the vehicle and the Sun being in opposite directions. The costate and Hamiltonian discontinuities at $t_d$ are given in Cartesian coordinates by

$$\begin{cases} \Delta\vec{p}_r(t_d) &= -\mu \dfrac{\partial \psi_S}{\partial \vec{r}} \\ \Delta\vec{p}_v(t_d) &= -\mu \dfrac{\partial \psi_S}{\partial \vec{v}} \\ \Delta p_m(t_d) &= -\mu \dfrac{\partial \psi_S}{\partial m} \\ \Delta H(t_d) &= \mu \dfrac{\partial \psi_S}{\partial t} \end{cases} \Rightarrow \begin{cases} \Delta p_x &= -\mu \sin \alpha_S \\ \Delta p_y &= \mu \cos \alpha_S \\ \Delta p_{vx} &= 0 \\ \Delta p_{vy} &= 0 \\ \Delta p_m &= 0 \\ \Delta H &= \mu \omega_S (x \cos \alpha_S + y \sin \alpha_S) \end{cases} \qquad (16)$$

where $\mu$ is an unknown multiplier associated to the interior point constraint Eq. (6).

To find the multiplier value, we write the Hamiltonian Eq. (12) just before and just after the discontinuity date. The state is continous. The only discontinous variables are the position costate $\vec{p}_r$ from Eq. (16) and the coefficient $\varepsilon$.

$$\begin{cases} H(t_d^-) &= \vec{p}_r(t_d^-) \cdot \vec{v} + \vec{p}_v \cdot \vec{g} + \varepsilon(t_d^-) T\Phi \\ H(t_d^+) &= \vec{p}_r(t_d^+) \cdot \vec{v} + \vec{p}_v \cdot \vec{g} + \varepsilon(t_d^+) T\Phi \end{cases} \qquad (17)$$

Substracting both equations yields a relationship between the Hamiltonian and the costate discontinuity.

$$\Delta H(t_d) = \Delta \vec{p}_r(t_d) \cdot \vec{v} + \Delta \varepsilon(t_d) T\Phi \qquad (18)$$

Using Eq. (16) to replace $\Delta H$ and $\Delta \vec{p}_r$, we get an explicit expression for the multiplier $\mu$.

$$\mu = \dfrac{T\Phi \Delta \varepsilon}{\dot{\psi}_S} \qquad (19)$$

where $\dot{\psi}_S = \dfrac{d\psi_S}{dt} = \dfrac{\partial \psi_S}{\partial t} + \dfrac{\partial \psi_S}{\partial \vec{r}} \cdot \vec{v}$ is the total derivative of the constraint function $\psi_S$ and $\Delta \varepsilon = -1$ (respectively $+1$) at the shadow entrance (respectively exit).

With this expression Eq. (19) the costate discontinuities Eq. (16) can be directly accounted within the propagation of the state and costate equations, provided that the dates of shadow entrance and exit are properly detected within the numerical integration.

The optimal trajectory can thus be found by solving the following two boundary value problem (TPBVP).

TPBVP

Find the initial costates $\vec{p}_r(t_0), \vec{p}_v(t_0), p_m(t_0)$, the final time $t_f$ and the constraints multipliers $\nu_A$, $\nu_P$ such that the transversality conditions Eqs. (14,15) and the final constraints Eq. (8) are met.

For a planar transfer the problem is of size 8. The shooting method consists in solving this nonlinear system by a Newton-like method with a numerical integration of the state and costate differential equations from the initial date to the final date. A major issue of this approach lies in the high sensitivity to the initial costate guess. Without a careful initialization and an accurate numerical integration the Newton method is very likely to fail. This behaviour is especially marked for low thrust transfers which address very long propagation durations (several days or weeks or even months) and dynamics discontinuities (due to eclipses). Finding a satisfying costate guess for such low thrust transfers is therefore a major challenge which may discourage from applying the shooting method.

## 3. Solution Method

The solution method proposed aims at bypassing the shooting method issues, due to the costate guess and to the numerical accuracy required by the Newton method.

### 3.1 Costate guess

The first part of the solution method consists in generating a correct costate guess. For that purpose we use some past results about orbital transfers.

A first useful result is that for minimum-time low thrust transfers the propellant consumption and the total velocity impulse are nearly insensitive to the thrust level [30,46-48]. When the transfer time becomes large the following relationship is observed.

$$m_0 - m_f = \frac{T t_f}{v_e} \approx C^{te} , \quad \forall T \tag{20}$$

A second useful result is that for minimum-fuel high thrust transfers at constant optimized thrust level, there exists an analytical costate solution, with orthogonal position and velocity costates [49].

$$\vec{p}_r(t_0) = \omega_0 \begin{pmatrix} -\cos\theta_0 \\ \sin\theta_0 \end{pmatrix} , \quad \vec{p}_v(t_0) = \begin{pmatrix} \sin\theta_0 \\ \cos\theta_0 \end{pmatrix} , \quad p_m(t_0) = \frac{v_e}{m_0} \tag{21}$$

In these formulae, $\theta_0$ is the initial pitch angle (between the horizontal and the thrust direction). It depends on the initial conditions (radius vector $r_0$, velocity modulus $v_0$, flight path angle $\gamma_0$) through the implicit equation

$$\sin(\theta_0 - \gamma_0) = \pm \frac{v_{c0}}{v_0} \frac{\sin\theta_0}{\sqrt{1 - 3\sin^2\theta_0}} \quad \text{with} \quad v_{c0} = \sqrt{\frac{\mu}{r_0}} \tag{22}$$

$\omega_0$ is the thrust initial rotation rate given by

$$\omega_0 = \sqrt{\frac{\mu}{r_0^3}(1 - 3\sin^2\theta_0)} \qquad (23)$$

The final mass for a low thrust transfer being nearly insensitive to the thrust level, it can be hoped that the above costate guess can be used successfully for low thrust transfers. We therefore consider the following costate guess parameterized by two angles $\theta_v$ and $\theta_n$

$$\vec{p}_r(t_0) = \omega_v \begin{pmatrix} -\cos\theta_n \\ \sin\theta_n \end{pmatrix} \;,\; \vec{p}_v(t_0) = \begin{pmatrix} \sin\theta_v \\ \cos\theta_v \end{pmatrix} \;,\; p_m(t_0) = \frac{v_e}{m_0} \qquad (24)$$

$$\omega_v = \sqrt{\frac{\mu}{r_0^3}(1 - 3\sin^2\theta_v)} \qquad (25)$$

For a low thrust transfer starting at the perigee, it is expected that the angle $\theta_v$ takes a small value (corresponding to a thrust aligned with the initial velocity) and that the angle $\theta_n$ remains approximatively equal to $\theta_v$. By this way the search space is reduced to a narrow interval of a few degrees around zero for both angles.

The angles $\theta_v$ and $\theta_n$ define the initial costate and therefore the complete command law along the trajectory using Eq. (11). The optimal control problem Eq. (8) is recast as a nonlinear programming problem (NLP) with 3 unknowns ($\theta_v$, $\theta_n$, $t_f$) and 2 constraints ($h_A$, $h_P$)

$$\min_{\theta_v, \theta_n, t_f} J = t_f \quad \text{s.t.} \quad \begin{cases} \psi_A[\vec{r}(t_f), \vec{v}(t_f)] = h_A \\ \psi_P[\vec{r}(t_f), \vec{v}(t_f)] = h_P \end{cases} \qquad (26)$$

### 3.2 Optimization method

The second part of the solution method consists in increasing the robustness to numerical inaccuracies. The propagation of the state and costate differential equations on long durations is prone to numerical errors. This inaccuracy is increased by the discontinuities due to eclipses. The numerous transitions at the shadow region bounds must be very accurately detected to ensure the smoothness of the integration result.
Noisy gradient assessments have very adverse effects for optimization algorithms [50] and more specifically when applying a shooting method. In order to cope with these numerical issues, we turn to a derivative free optimization method. Such methods can be very efficient on "noisy" optimization problems suffering from numerical inaccuracy or poor smoothness, but they are generally not suited to constraint handling [51,52].

In our specific case of a low thrust orbital transfer, it is possible to reduce the NLP problem Eq. (26) to a pure minimization problem without constraints. In most practical cases, the transfer aims at raising the vehicle up to a high altitude orbit such as the geostationary orbit (GEO). For such transfers the perigee altitude is monotonously increasing along the trajectory and the transfer is achieved once the perigee reaches the targeted altitude $h_P$. This property can be exploited to discard the final date from the problem unknowns. The motion integration is stopped as soon as the perigee reaches the desired value $h_P$. The cost function is the squared difference between the final and targeted apogees to minimize.

The NLP problem reduces thus to

$$\min_{\theta_v, \theta_n} \left(\psi_A[\vec{r}(t_f), \vec{v}(t_f)] - h_A\right)^2 \quad (27)$$

It is not mathematically guaranteed that this formulation actually yields the minimum-time trajectory. Indeed different sets of $\theta_v$ and $\theta_n$ values could exist that allow reaching the targeted orbit with different final times. In practice the narrow interval of search for these angles (with an initial thrust direction making a small angle with the velocity) and the trajectory plots give a good confidence in the result optimality.

The formulation Eq. (27) with 2 unknowns is suited to any derivative free optimization method [53]. A particularly attractive algorithm for this small size problem is DIRECT [54-56]. DIRECT explores the search space in a deterministic way and ensures finding the global minimum under reasonable regularity assumptions (Lipschitzian cost function). Opposite to metaheuristics like genetic algorithm, particle swarm, simulated annealing, … that involve a stochastic part and therefore a risk of missing the global minimum, DIRECT is able to locate correctly the solution with a controlled number of function evaluations.

With a derivative free algorithm it is not necessary to pay a special care on the numerical integration accuracy. For the present application, the state and costate are expressed in Cartesian coordinates. This choice is not numerically the best for a low thrust orbital transfer, but it yields very simple equations for the state equations, the costate equations and the eclipse constraint formulation. The eighth order Runge-Kutta method DOP853 is used for the propagation with large time steps and with a large tolerance on the integration error. The eclipses are detected within the propagation and the time step is adapted by dichotomy to pass accurately at the eclipse transition date and apply the explicit costate discontinuity given by Eqs. (16,19). These tunings make each simulation quite fast. The convergence of DIRECT requires a few thousands simulations and it is achieved in a few minutes in the longest cases (transfer durations of several months).

## 4. Application

The solution method is illustrated on a coplanar low thrust transfer from a geostationary transfer orbit (GTO) to the geostationary orbit (GEO) with different thrust level values and different initial lightening conditions.

### 4.1 Illustrative Example

The vehicle initial gross mass is 1000 kg, the engine thrust level is 1 N with a specific impulse of 1500 s corresponding to an exhaust velocity of 14710 m/s.

The initial orbit is a GTO with apogee and perigee altitudes respectively at 36000 km and 500 km. The x axis of the reference Galilean frame is defined by the initial perigee. The initial anomaly is thus equal to 0 deg. The perigee initial local time is 0 h and the Sun direction moves at the angular rate of 0.986 deg/day.

The target is a circular GEO at 36000 km. The transfer is planar. The Earth equatorial radius is $R_E = 6378137$ m and the gravitational constant is $\mu = 3.986005 \cdot 10^{14}$ m$^3$/s$^2$.

The derivative free algorithm DIRECT is used to find the 2 angles $\theta_v$ and $\theta_n$ orientating the initial costate. The search interval is [−10 deg ; +10 deg] for both angles. With a time step of 1000 s and accuracy of 10 s for the eclipse detection, each simulation takes on average 0.02 s on a standard Linux computer using the dynamics equations in Cartesian coordinates. The convergence is achieved within 1 min and it requires about 3000 function calls.

The Figure 2 presents the minimum-time trajectories found with a 1 N engine assuming initial perigee local times of respectively 0 h (Sun above apogee), 6 h and 12 h (Sun above perigee). The shadow region (null thrust) is represented by black bold lines. It rotates simultaneously with the apparent Sun motion at a 0.986 deg/day rate. The thrust direction is indicated by red arrows during the first revolution and green arrows during the last one. It can be seen on the angle of attack plot that the thrust is mostly directed along the velocity during the first part of the transfer making the apogee increase beyond 36000 km. During the second part of the transfer the thrust direction oscillates at the orbital period, along the velocity at the apogee, and opposite to the velocity at the perigee to bring the apogee back to 36000 km. The eclipse duration ranges from 21 h when the apogee is lightened to 70 h when the apogee is in the Earth shadow. This last configuration is the worst since a large time is spent without thrusting and moreover the optimal thrusting location (the apogee) is forbidden.

The values of the final time, of the total eclipse duration and of the optimized angles $\theta_v$ and $\theta_n$ are given on the right side of Figure 2. It must be recalled that these two angles defines completely the initial costate vector, so that only sub-optimal solutions can be found. The angles values are close to zero as expected from the analytical costate solution derived for high thrust level transfers.

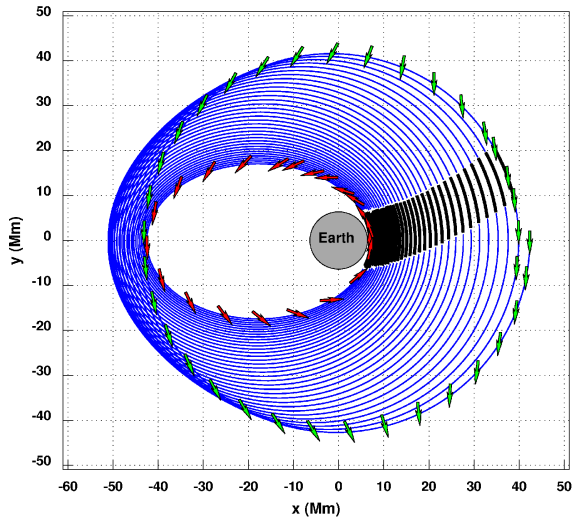
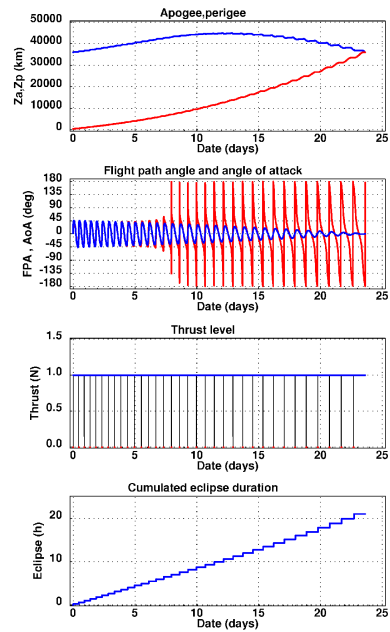

Figure 2 :
Transfers with 1N

| Perigee local time | 0 h |
|---|---|
| Transfer duration | 23.63 day |
| Eclipse duration | 21.1 h |
| $\theta_v$ | 0.0004 deg |
| $\theta_n$ | −0.0102 deg |

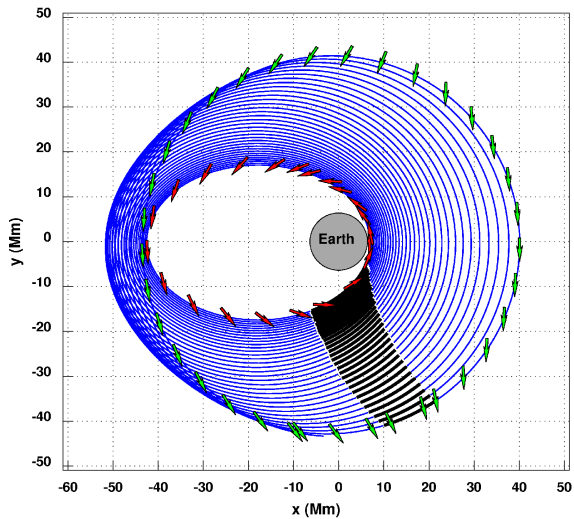
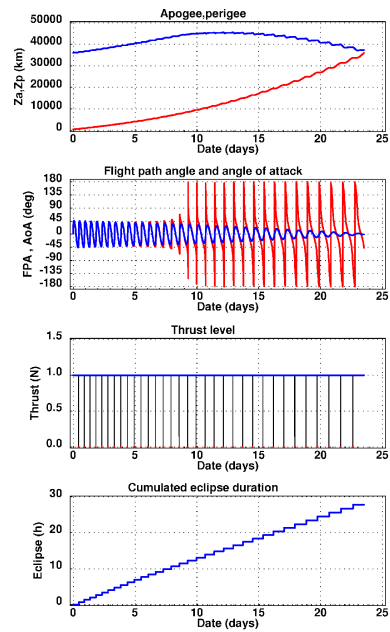

| Perigee local time | 6 h |
|---|---|
| Transfer duration | 23.50 day |
| Eclipse duration | 27.9 h |
| $\theta_v$ | 1.399 deg |
| $\theta_n$ | −2.477 deg |

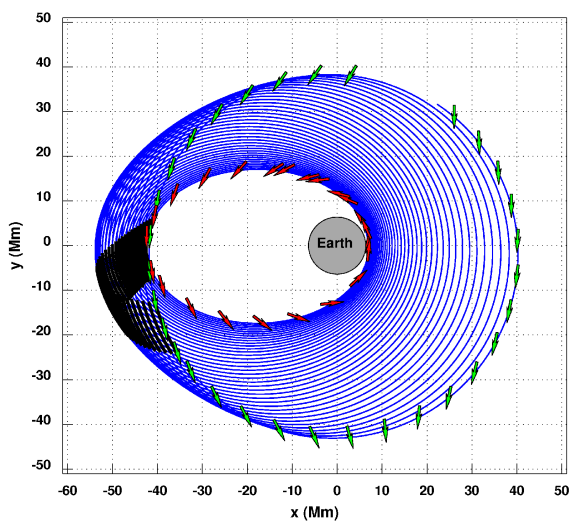
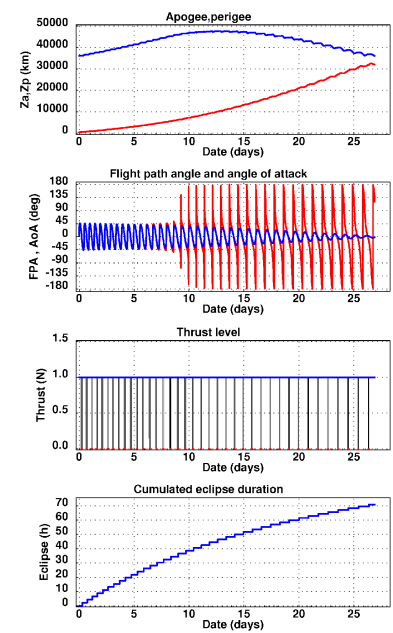

| Perigee local time | 12 h |
|---|---|
| Transfer duration | 27.01 day |
| Eclipse duration | 69.6 h |
| $\theta_v$ | 3.364 deg |
| $\theta_n$ | −5.949 deg |

### 4.2 Sensitivities

The above example is solved again for thrust levels of 10 N and 0.1 N, and for initial perigee local times ranging from 0 h to 12 h. The Table 1 summarizes the numerical results obtained. The transfer without eclipse is also assessed and compared to results given by the CNES software MIPELEC [57,58].

| Thrust level (N)       | 10    |       |       |
|------------------------|-------|-------|-------|
| Perigee local time (h) | 0 h   | 6 h   | 12 h  |
| Final mass (kg)        | 863,7 | 867,0 | 858,2 |
| Final time (day)       | 2,42  | 2,39  | 2,72  |
| Velocity impulse (m/s) | 2155,7| 2100,1| 2249,5|
| Number of revolutions  | 3     | 3     | 4     |
| Eclipse duration (h)   | 2,3   | 2,9   | 7,4   |

| Thrust level (N)       | 1     |       |       |       |       |           |         |
|------------------------|-------|-------|-------|-------|-------|-----------|---------|
| Perigee local time (h) | 0 h   | 3 h   | 6 h   | 9 h   | 12 h  | No eclipse| MIPELEC |
| Final mass (kg)        | 866,4 | 868,1 | 868,8 | 865,4 | 858,4 | 864,0     | 864,3   |
| Final time (day)       | 23,63 | 23,38 | 23,50 | 24,78 | 27,01 | 23,15     | 23,10   |
| Velocity impulse (m/s) | 2109,7| 2081,6| 2069,4| 2126,7| 2246,0| 2150,1    | 2145.2  |
| Number of revolutions  | 33    | 33    | 33    | 35    | 37    | 32        | 32      |
| Eclipse duration (h)   | 21,1  | 22,1  | 27,9  | 44,7  | 69,6  |           |         |

| Thrust level (N)       | 0,1    |        |        |        |        |           |         |
|------------------------|--------|--------|--------|--------|--------|-----------|---------|
| Perigee local time (h) | 0 h    | 3 h    | 6 h    | 9 h    | 12 h   | No eclipse| MIPELEC |
| Final mass (kg)        | 867,2  | 867,1  | 868,8  | 869,4  | 865,7  | 864,1     | 864,3   |
| Final time (day)       | 239,40 | 236,95 | 233,20 | 233,89 | 245,90 | 231,33    | 231,03  |
| Velocity impulse (m/s) | 2096,2 | 2097,4 | 2068,1 | 2059,3 | 2121,3 | 2148,2    | 2145.2  |
| Number of revolutions  | 341    | 328    | 328    | 330    | 354    | 319       | 319     |
| Eclipse duration (h)   | 317,7  | 257,4  | 238,2  | 274,9  | 414,1  |           |         |

Table 1 : GTO-GEO transfer with various thrust levels and initial perigee times

An apogee in the Earth shadow (perigee time = 12 h) gives the worst results since the most favorable part of the orbit is not available for thrusting. Perigee times between 0 h and 6 h yield near performances both in terms of transfer duration and final mass. These two quantities are not directly linked since the eclipse location both influences the transfer duration and the thrust efficiency wrt the orbit evolution. The best performances are obtained for a Sun location orthogonal to the initial semi major axis (perigee time = 6 h). In that configuration the shadow area does not affect the apsides which are the most efficient thrusting places. This results simultaneously in a reduction of the transfer duration and in an increase of the final mass.

On the other hand a good agreement is found with MIPELEC results which are obtained by an averaging method. MIPELEC yields the same final mass whatever the thrust level and the final time is deduced from the propellant consumption divided by the propulsive mass flow-rate.

Comparing the transfer duration with and without eclipse it is also observed that the rule of the thumb "difference = total eclipse duration" would be rather pessimistic. In most cases the duration penalty due to eclipses is largely lower than the time spent in eclipse.

The Figure 3 presents the trajectories obtained for thrust levels of 10 N and 0.1 N. For the 0.1 N case, the shadow region is not colored in order to better distinguish it.

The method is also applied to a LEO to GEO transfer with an initial circular orbit at 500 km, considering the same vehicle assumptions (mass = 1 t, thrust level = 1 N, Isp = 1500 s). The results are presented in the Table 2. The comparison with MIPELEC is still correct.

| Thrust level (N)       | 1      |            |         |
|------------------------|--------|------------|---------|
| Perigee local time (h) | 0 h    | No eclipse | MIPELEC |
| Final mass (kg)        | 725,5  | 734,2      | 734,2   |
| Final time (day)       | 58,30  | 45,26      | 45,26   |
| Velocity impulse (m/s) | 4720,3 | 4545,0     | 4545,0  |
| Number revolutions     | 398    | 297        | 297     |
| Eclipse duration (h)   | 277,9  |            |         |

Table 2 : LEO-GEO transfer with 1 N

### 4.3 Extensions

The derivative free approach combined with the analytical costate guess proves satisfying for coplanar low-thrust transfers. The solution is found from scratch in a few minutes at most without requiring specific guesses or tunings from one case to another. The extension to three dimensional transfers is not straightforward, and the following questions should be addressed :

- Is it still possible to build an adequate analytical costate guess to reduce the number of unknowns ?
- Is it still possible to derive explicit costate discontinuities at the eclipse bounds in order to apply them directly within the trajectory integration ?
- How to define the cost function of the unconstrained problem since the apogee value and the inclination have to be simultaneously targeted ?
- Does the DIRECT algorithm still behave well when the number of unknowns and the search space increase ?

These questions are currently under study with the aim of extending the solution method to three dimensional low-thrust transfers.

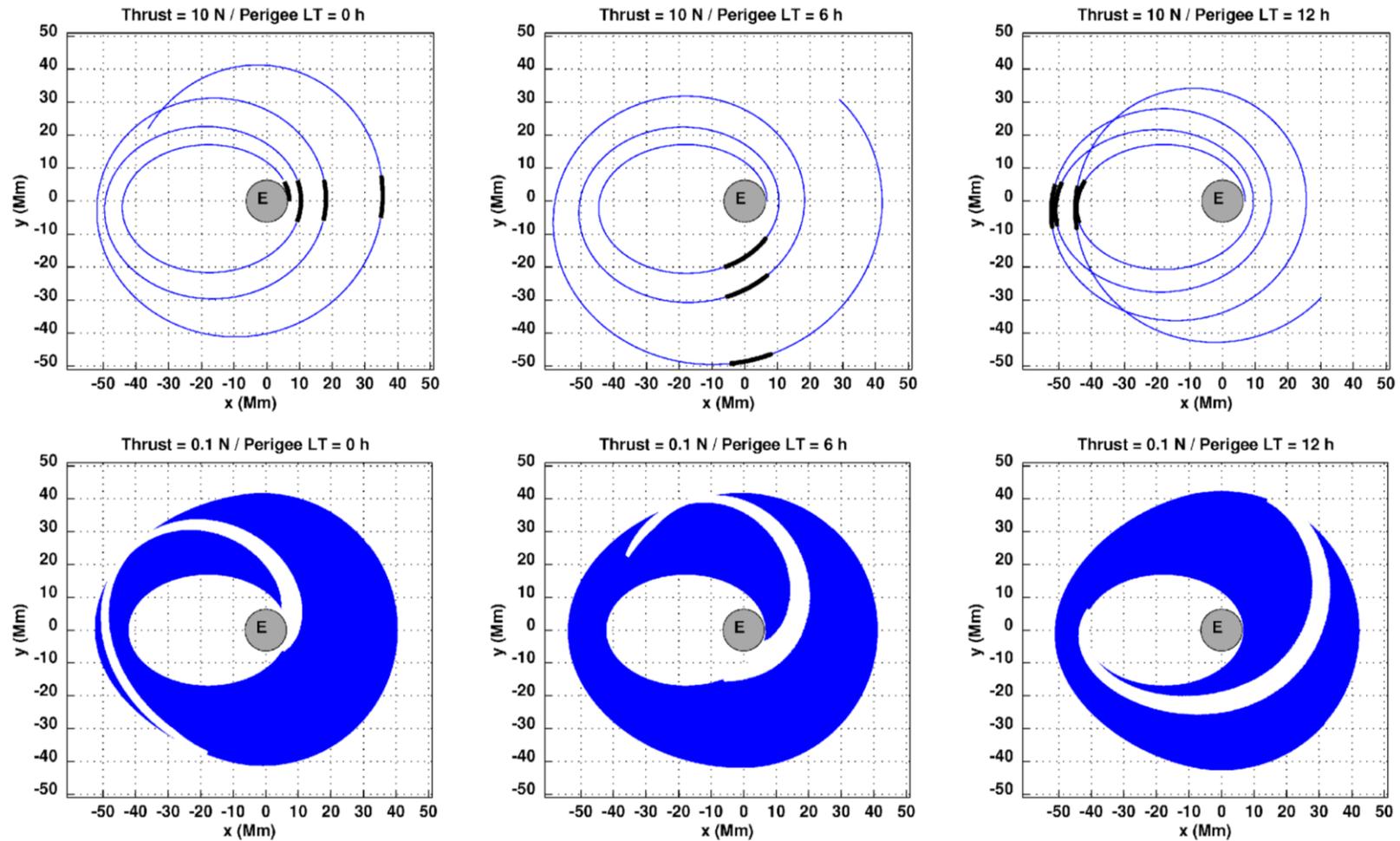

Figure 3 : Transfers with 10 N and 0.1 N

# 5. Conclusion

The minimum-time low thrust transfer with eclipses has been investigated in the coplanar case. The trajectory is simulated aligning the thrust direction with the velocity costate, as given by the Pontryaguin Maximum Principle. An explicit expression for the costate discontinuities at the eclipse entrances and exits is derived. These discontinuities are directly accounted within the numerical integration of the trajectory. An analytical costate solution obtained for high thrust transfers suggests parameterizing the initial costate vector by two angles. This assumption proves satisfying and it allows reducing the shooting problem to an unconstrained nonlinear problem with two unknowns angles taking near zero values. Applying the derivative free algorithm DIRECT avoids dealing with numerical accuracy issues for the trajectory integration and gives some guarantee of locating the global minimum in the narrow search space.. The method has been experimented on transfers towards the geostationnary orbit starting from different initial orbits, using different thrust levels from large to very low, and considering different initial lightening conditions. In all cases a near optimal solution is found from scratch in a few minutes at most without requiring specific guesses or tunings from one case to another. The extension of the method to three dimensional transfers with inclination changes requires significant improvements of the solution method and it is currently under study.

**References**

[1] Betts, J. T., Survey of numerical methods for trajectory optimization, Journal of Guidance, Control, and Dynamics, vol. 21, no. 2, pp. 193–207, 1998.
[2] Rao, A., A survey of numerical methods for optimal control, Advances in the Astronautical Sciences, vol. 135, no. 1, pp. 497–528, 2009.
[3] Kim, M., Continuous Low-Thrust Trajectory Optimisation : Techniques and Applications, PhD thesis, Virginia Polytechnic Institute and State University, 2005.
[4] Hargraves, C. R., Paris, S. W., Direct trajectory optimization using nonlinear programming and collocation, Journal of Guidance, Control, and Dynamics, vol. 0, no. 4, pp. 338–342, 1987.
[5] Enright, P. J., Conway, B. A., Discrete approximations to optimal trajectories using direct transcription and nonlinear programming, Journal of Guidance, Control, and Dynamics, vol. 15, no. 4, pp. 994–1002, 1992.
[6] Ross, I. M., Fahroo, F., Pseudospectral knotting methods for solving optimal control problems, Journal of Guidance, Control, and Dynamics, vol. 27, no. 3, pp. 397–405, 2004.
[7] Paris, S. W., Riehl, J. P., Sjauw, W. K., Enhanced procedures for direct trajectory optimization using nonlinear programming and implicit integration, in Proceedings of the AIAA/AAS Astrodynamics Specialist Conference and Exhibit, pp. 21–24, August 2006.
[8] Dileep M V, Vishnu G Nair, Prahalad K R, Surekha Kamath, Trajectory Optimization of Launch Vehicles Using Steepest Descent Method – A Novel Approach Ascent Optimization for a Heavy Space Launcher, Int. Journal of Engineering Research and Applications, ISSN : 2248-9622, Vol. 4, Issue 1( Version 1), January 2014, pp.116-121.
[9] Betts, J.T., Practical methods for optimal control and estimation using nonlinear programming, Siam, 2010.
[10] Topputo, F., Zhang, C., Survey of Direct Transcription for Low-Thrust Space Trajectory Optimization with Applications, Abstract and Applied Analysis, Volume 2014, Article ID 851720, DOI: 10.1155/2014/851720.
[11] Betts, J. T., Very low thrust trajectory optimization using a direct SQP method, Journal of Computational and applied Mathematics, vol. 120, issues1-2, pp. 27–40, Aug. 2000.
[12] Betts, J.T., Optimal Low Thrust Orbit Transfer with Eclipsing, 14 January 2014.
[13] Lawden, D.F., Optimal transfer between coplanar elliptical orbits, Journal of Guidance Control and Dynamics, vol. 15, no. 3, pp. 788-791, 1991.

[14] Broucke, R.A., Prado, A.F.B.A., Optimal N-impulse transfer between coplanar orbits , AAS/AIAA Astrodynamics Specialist Conference, Victoria, Canada, 1994.
[15] Yam, C.H., Izzo, D., Biscani, F., Towards a High Fidelity Direct Transcription Method for Optimisation of Low-Thrust Trajectories, ESTEC, 2011.
[16] Landau, D.F., Longuski, J.M., Trajectories for Human Missions to Mars, Part 2 : Low Thrust Transfers, Journal of Spacecrafts and Rockets, Vol. 43, No. 5, September-October 2006.
[17] Pontryagin, L., Boltyanskii, V., Gramkrelidze, R., Mischenko, E., The mathematical theory of optimal processes, Wiley Interscience, 1962.
[18] Trélat, E., Contrôle optimal – Théorie et Applications, Vuibert, 2005.
[19] Conway, B.A., Spacecraft Trajectory Optimization, Cambridge University Press, 2010.
[20] Lawden, D.F. Optimal trajectories for space navigation, Butterworths Publishing Corporation, 1963.
[21] Leitmann, G., Optimization techniques with applications to Aerospace System, Academic Press New-York, 1962.
[22] Marec, J.P., Optimal space trajectories, Elsevier (1979)
[23] Bonnard, B., Caillau, J.-B., Trélat, E., Geometric optimal control of elliptic Keplerian orbits, Discrete Cont. Dyn. Syst. 5 (2005), no. 4, 929--956.
[24] Caillau, J.-B., Noailles, J., Coplanar control of a satellite around the Earth, ESAIM Cont. Optim. Calc. Var. 6 (2001), 239--258.
[25] P. Martinon, J. Gergaud, Using switching detection and variational equations for the shooting method, Optimal Control and Applications and Methods 28, no. 2 (2007), 95--116.
[26] Bonnans, F., Martinon, P., Trélat, E., Singular arcs in the generalized Goddard's problem , Journal of Optimization Theory and Applications, Vol 139 (2), pp 439-461, 2008.
[27] Ponssard, P., Graichen, K., Petit, N., Laurent-Varin, J. Ascent Optimization for a Heavy Space Launcher, Proceedings of the European Control Conference 2009, Budapest, Hungary, August 23-26, 2009.
[28] Augros, P., Delage, R., Perrot, L., Computation of optimal coplanar orbit transfers, AIAA 1999.
[29] Naidu, D.S., Hibey, J.L., Charalambous, C., Fuel-Optimal Trajectories for Aeroassisted Coplanar Orbital Transfer Problem, IEEE Transactions on Aerospace and Electronic Systems 26 no.2 (March 1990) 374-381.
[30] Gergaud, J., Haberkorn, T., Orbital transfer : some links between the low-thrust and the impulse cases, Acta Astronautica, 60, no. 6-9 (2007), 649-657.
[31] Lee, D., Bang, H., Efficient Initial Costates Estimation for Optimal Spiral Orbit Transfer Trajectories Design, Journal of Guidance, Control and Dynamics, Vol. 32, No. 6, November–December 2009
[32] Pifko, S.L., Zorn, A.H., West, M., Geometric Interpretation of Adjoint Equations in Optimal Low Thrust Space Flight, AIAA/AAS Astrodynamics Specialist Conference and Exhibit, 18-21 August 2008, Hawaii.
[33] Oberle, H.J., Taubert, K. , Existence and multiple solutions of the minimum-fuel orbit transfer problem, Journal of Optimization Theory Appl. 95 (1997) 243-262.
[34] Graham, K.F., Rao, A.V., Minimum-Time Trajectory Optimization of Low Thrust Earth-Orbit Transfers with Eclipsing, AAS/AIAA Space Flight Mechanics Meeting, Williamsburg, Virginia, 11-15 January 2015.
[35] Jiang, F., Baoyin, H., Li, J., Practical techniques for low-thrust trajectory optimization with homotopic approach, Journal of Guidance, Control, and Dynamics, vol. 35, no. 1, pp. 245–258, 2012.
[36] Gergaud, J., Haberkorn, T., Martinon, P., Low thrust minimum fuel orbital transfer : an homotopic approach, Journal of Guidance Control and Dynamics, 27, 6 (2004), 1046-1060.
[37] Thorne, J.D., Hall, C.D., Minimum-Time Continuous Thrust Orbit Transfers, Journal of the Astronautical Sciences 47 (1997), no. 4, 411-432.
[38] Dargent, T., Martinot, V., An integrated tool for low thrust optimal control orbit transfers in interplanetary trajectories, Proceedings of the 18[th] International Symposium on Space Flight Dynamics (ESA SP-548), 2004, 6 pages.
[39] Pontani, M., Conway, B.A., Particle Swarm Optimization Applied to Space Trajectories, Journal of Guidance, Control, and Dynamics, Vol. 33, No. 5 (2010), pp. 1429-1441, DOI: 10.2514/1.48475.
[40] Steffens, M.J., A combined global and local methodology for launch vehicle trajectory design-space exploration and optimization, Thesis Georgia Institute of technology, April 2014
[41] Zuiani, F., Multi-objective optimization of low-thrust trajectories. PhD thesis, University of Glasgow, 2015.
[42] Bérend, N., Bonnans, F., Haddou, M., Laurent-Varin, J., Talbot, C., An Interior-Point Approach to Trajectory Optimization, Journal of Guidance, Control and Dynamics 30 (2007), no. 5, 1228--1238.
[43] Bourgeois, E., Bokanowski, O., Zidani, H., Désilles A., Optimization of the launcher ascent trajectory leading to the global optimum without any initialization: the breakthrough of the HJB approach, 6[th] European Conference for Aeronautics and Space Sciences (EUCASS) 29 June-2 July 2015.
[44] Bryson, A.E., Ho, Y.C. Applied optimal control, Hemisphere Publishing Corporation, 1975.
[45] Hull, D.G., Optimal control theory for applications, Springer, 2003.
[46] Vallado, D. A., Fundamentals of Astrodynamics and Applications, Third Edition, Microcosm Press, El Segundo, California, 2007.

[47] Bombrun, A., Pomet, J.-B., Asymptotic Behavior of Time Optimal Orbital Transfer for Low Thrust 2-Body Control System, Discrete and Continuous Dynamical Systems Supplement 2007, pp 122-129.
[48] Thorne, J.D., Minimum-Time Constant-Thrust Orbit Transfers with Noncircular Boundary Conditions, 2008.
[49] Cerf, M., Optimal Thrust Level for Orbit Insertion, Arxiv preprint, 2016.
[50] Vugrin, K.E., On the Effect of Numerical Noise in Simulation-Based Optimization, Faculty of the Virginia Polytechnic Institute and State University, Thesis, March 2003.
[51] Di Pillo, G., Lucidi, S., Rinaldi, F., A Derivative-Free Algorithm for Constrained Global Optimization Based on Exact Penalty Functions, Journal of Optimization Theory and Applications, DOI 10.1007/s10957-013-0487-1, November 2013.
[52] Parsopoulos, K.E., Vrahatis, M.N., Particle Swarm Optimization Method for Constrained Optimization Problems, Frontiers in Artificial Intelligence and Applications, January 2002.
[53] Paluszek, M.A., Thomas, S.J., Trajectory Optimization Using Global Methods, IEPC05-173, 29$^{th}$ International Electric Propulsion Conference, Princeton University, 31 October – 4 November, 2005.
[54] Perttunen, C.D., Jones, D.R., Stuckman, B.E., Lipschitzian optimization without the lipschitz constant, Journal of Optimization Theory and Application, 79(1):157-181, October 1993.
[55] Gablonsky, J.M., Direct version 2.0 userguide. Technical Report, CRSC-TR01-08, Center for Research in Scientific Computation, North Carolina State University, April 2001.
[56] Liuzzi, G., Lucidi, S., Piccialli, V., Exploiting Derivative-Free Local Searches In Direct-Type Algorithms for Global Optimization, Istituto di Analisi dei Sistemi ed Informatica "A. Ruberti", Rome, Italy, ISSN: 1128–3378, 2014.
[57] Geffroy, S., Généralisation des techniques de moyennation en contrôle optimal - Application aux problèmes de transfert et rendez-vous orbitaux à poussée faible, Thèse de doctorat, Institut National Polytechnique de Toulouse, 1997.
[58] Geffroy, S., Epenoy, R., Optimal low-thrust transfers with constraints - Generalization of averaging techniques, Acta Astronautica, Vol. 41, No. 3, 1997, pp. 133–149.

**Acronyms**

| | |
|---|---|
| OCP | Optimal Control Problem |
| PMP | Pontryaguin Maximum Principle |
| TPBVP | Two Point Boundary Value Problem |
| NLP | Non Linear Programming |
| DIRECT | Dividing RECTangles |
| GEO | Geostationary Orbit |
| GTO | Geostationary Transfer Orbit |
| LEO | Low Earth Orbit |